\def\timestamp{%
Time-stamp: <alan-arxiv.tex: Wednesday 30-11-2016 at 11:22:10 (cet)>}
\def\stripname Time-stamp: <#1: #2 #3 at #4 #5>{#3 -  #4 #5}
\edef\filedate{\expandafter\stripname\timestamp}
\documentclass[a4paper]{amsart}
\usepackage{amsrefs}

\DeclareMathSymbol\HH0{AMSb}{`H}
\DeclareMathSymbol\N0{AMSb}{`N}
\DeclareMathSymbol\Q0{AMSb}{`Q}
\DeclareMathSymbol\R0{AMSb}{`R}

\newcommand\betaN{\beta\N}

\newcommand\betaR{\beta\R}

\newcommand\Hstar{\HH^*}
\newcommand\Nstar{\N^*}
\newcommand\Qstar{\Q^*}
\newcommand\Rstar{\R^*}
\newcommand\bee{\mathfrak{b}}
\newcommand\cee{\mathfrak{c}}
\newcommand\ess{\mathfrak{s}}

\let\sbs\setminus
\let\bold\mathbf
\let\Bbb\mathbb

\newcommand\cf{\operatorname{cf}}
\newcommand\cl{\operatorname{cl}}
\newcommand\clb{\cl_\beta}
\newcommand\cls{\cl_s}

\newcommand\axiom{\mathsf}
\newcommand\CH{\axiom{CH}}
\newcommand\PFA{\axiom{PFA}}
\newcommand\ZFC{\axiom{ZFC}}

\newcommand\calA{\mathcal{A}}
\newcommand\calD{\mathcal{D}}

\usepackage{graphicx}

\begin{document}

\begin{center}
  \includegraphics[width=5cm]{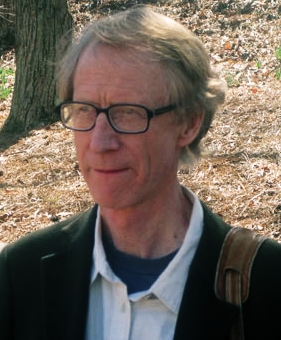}
\end{center}

\title{Alan Dow}

\author[K. P. Hart]{Klaas Pieter Hart}
\address[K. P. Hart]%
        {Faculty of Electrical Engineering, Mathematics and Computer Science\\
          TU~Delft\\
          Postbus 5031\\
          2600~GA Delft\\
          The Netherlands}
\email{k.p.hart@tudelft.nl}
\urladdr{http://fa.its.tudelft.nl/\~{}hart}

\author[J. van Mill]{Jan van Mill}
\address[J. van Mill]{KdV Institute for Mathematics\\
                      University of Amsterdam\\
                      Science Park 904\\
                      P.O. Box 94248\\
                      1090~GE Amsterdam\\
                      The Netherlands}
\email{j.vanmill@uva.nl}
\urladdr{https://staff.fnwi.uva.nl/j.vanmill}
\address[J. van Mill]%
        {Faculty of Electrical Engineering, Mathematics and Computer Science\\
         TU~Delft\\ Postbus 5031\\ 2600~GA Delft\\ The Netherlands}

\subjclass{03E, 54}
\date{\filedate}

\maketitle

Alan Dow celebrated his 60th birthday on the 5th of December 2014 by flying
to Ithaca, NY, to undergo a conference in honour of that birthday.
This article is an attempt to present an overview of Alan's many contributions
to General and Set-theoretic Topology.

\section{Elementarity}

One of the many gifts from Alan to Set-theoretic Topology is the use
of elementarity.
For a while this was even known as ``Dow's method of elementary submodels''.
But Alan would be, was, and still is the first to protest that
the L\"owenheim-Skolem theorem predates him by a few decades.

We have for the longest time been familiar with recursive constructions
where often beforehand a sequence of situations\slash sets is set up
and during the construction witnesses to bad things will be eliminated
or witnesses to good things will be embraced.
In the end we consider such a situation and realize that it was basically
dealt with during the construction.
A very good example is the Pol-Shapirovski\u\i{} proof of Arhangel$'$ski\u\i's
theorem on the cardinality of compact first-countable spaces.

What set theorists realised was that one can reduce the length of such proofs
considerably by an application of the L\"owenheim-Skolem theorem to a model
of `enough set theory': its proof is the ultimate closing-off argument where
one deals with \emph{all possible} situations in one go (even ones that
will never occur in your problem at hand).
But, and this is where this method gets its power, you will certainly have
dealt with every eventuality related to your problem.
Basically what is left is to perform what would have been the final step of
your old recursive argument.
This requires some familiarity with first-order logic and model theory,
so that you know how far you can go with your arguments.
But the time spent learning that will pay itself back handsomely in time
saved later.

A good place to start learning this is Alan's first introduction,
\cite{MR1031969}, which has an elementary proof of Arhangel$'$ski\u\i's
theorem that one should put next to the Pol-Shapirovski\u\i{} argument
to see difference between the `standard' and the elementarity mindset.
A later survey, \cite{MR1342520}, gives more applications of the latter.

\section{Remote points}

If $X$ is completely regular then it is dense in $\beta X$, so every point
of~$X^*$ lies close to~$X$; however some points lie closer to~$X$ than others.
One can formulate degrees of closeness by stipulating that the point belongs
to the closure of a topologically small subset of~$X$.
Thus, for example, $p\in X^*$ is \emph{near} if $p\in\clb D$ for some
closed and discrete subset of~$X$; other variants can be obtained by
using relatively discrete subsets, scattered subsets, and nowhere dense
sets.
The negation of the last notion has proved to be very fruitful:
call $p\in X^*$ a \emph{remote point of~$X$} if $p\notin\clb A$
for all nowhere dense subsets of~$X$.
Remote points were introduced by Fine and Gillman in 1962 who proved that the
Continuum Hypothesis (abbreviated $\CH$) implies that the real line~$\R$ has a
remote point.
Actually, their proof applies to every separable and non-pseudocompact space.
Around~1980 van~Douwen, and independently Chae and Smith, proved in~$\ZFC$
that every non-pseudocompact space with countable $\pi$-weight has remote
points.
That there are spaces without remote points, was demonstrated by
van~Douwen and van~Mill.
In 1982, Alan took over the research on remote points completely,
leaving absolutely nothing for his competitors
(see
\cites{MR694543,MR702616,MR774511,MR969061,MR678677,MR983872,MR1971307,
       MR2495242,MR3414887}).
He became the world's expert on remote points.
He substantially improved the results of van~Douwen, and Chae and~Smith
by showing that every non-pseudocompact ccc space of $\pi$-weight at
most~$\omega_1$ has a remote point,
and that under~$\CH$ the bound~$\omega_1$ is not optimal.
The fruits of remote points are manifold.
The points themselves were used in `honest' proofs of non-homogeneity of
certain \v{C}ech-Stone remainders: for example $\Qstar$~is extremally
disconnected at each remote point but not at other points.
The techniques developed and used for their construction have found many
applications too.
Alan's proof that $\omega\times 2^\kappa$ has remote points gave new insight
in the structure of the partial order that adds Cohen reals: a remote point,
seen as a clopen filter on $\omega\times2^\kappa$, takes big bites out of
dense open sets and these bites combine to form approximations of generic
filters, called \emph{enDowments} by some.
These enDowments were crucial in a Cohen-real proof of the consistency of
the normal Moore space conjecture, \cites{MR1075372,MR1080345}.
We must also mention that Alan showed that the result of Fine and Gillman
needs extra assumptions: in the side-by-side Sacks model there is a
non-pseudocompact separable space without remote points, \cite{MR983872}.

\section{$\betaN$ and $\betaR$}

The \v{C}ech-Stone compactification $\betaN$ of the discrete space of
natural numbers, $\N$, is Alan's favorite space, just like the space of
rational numbers, $\Q$, was Eric van~Douwen's favorite space.
Alan's work concentrates on the construction of special ultrafilters on~$\N$
such as weak $P$-points, cozero-accessible points, (bow)tie-points,
certain (finite-to-one) maps defined on $\betaN$ or
$\Nstar=\betaN\setminus\N$, and closed subsets of $\Nstar$.
Although Alan can be thought of as a set theorist, his thinking is topological:
he is interested in those set theoretical aspects of $\betaN$ that shed light
on topological questions.

He demonstrated his $\betaN$-talents already quite early in his career.
In~\cite{MR1277880} he gave an overview of certain aspects of $\betaN$
that is still valuable today.
His proof presented there that there are $2^\cee$ Rudin-Keisler incomparable
ultrafilters on~$\N$, is the best one around. 

The space $\betaN$ contains many copies of itself.
If $D$~is any countable discrete subset of~$\betaN$, then its closure is a
(topological) copy of~$\betaN$.
Hence the space~$\Nstar$ contains many topological copies of itself as well.
This prompted Eric van~Douwen to ask whether there is a copy of~$\Nstar$
in~$\Nstar$ that is not of the form $\cl D\setminus D$ for some countable
discrete subset~$D$ of~$\Nstar$.
Such a copy of~$\Nstar$ is called \emph{nontrivial}.
Under various assumptions, nontrivial copies of~$\Nstar$ exist.
But what about such copies in~$\ZFC$?
This longstanding open problem, which was thought to be beyond reach by
experts, was solved in the affirmative by Alan in~\cite{MR3209343}.

In~\cite{MR1752102}, he showed with Hart that the Open Coloring Axiom implies
that the Stone space of the measure algebra is not a continuous image
of~$\Nstar$, which contradicts an earlier result of Frankiewicz and Gutek.
Alan's work related to van~Douwen
spaces~\cites{MR2365321,MR2361681,MR2282377,MR2162384} gives answers to
very natural problems on spaces that are continuous images of~$\betaN$
or~$\Nstar$ under finite-to-one mappings.

One of Alan's favourite results about $\betaN$ is the existence
of a tree $\pi$-base for~$\Nstar$, so naturally he would investigate
the possible structure of these, \cite{MR1020980}, as well as employ
them in the study of the absolute of $\Nstar$, \cites{MR759135,MR1619290}.
The latter paper ties in with another \v{C}ech-Stone compactification
that has been the beneficiary of Alan's interest: $\betaR$.
As with $\betaN$, it is mostly $\Rstar$, or rather one of its
halves, $\Hstar$, whose structure we would like to clear
up --- $\HH=[0,\infty)$. 

One thing to do is look for parallels between $\Nstar$ and $\Hstar$:
take a known result on~$\Nstar$ and reformulate it to take into account
that $\Hstar$~is a continuum.
Sometimes this works, as in~\cite{MR1707489}:
there is a complete parallel version of Parovi\v{c}enko's theorem:
every continuum of weight~$\aleph_1$ or less is a continuous image
of~$\Hstar$ and thus, under~$\CH$, the continua
of weight~$\cee$ or less are exactly the continuous images of~$\Hstar$.
The parallel also extends in the negative direction: many examples
of non-images of~$\Nstar$ have a connected counterpart.

Sometimes the parallel breaks down: every separable compact space is clearly
a continuous image of~$\Nstar$ but there is a separable continuum that,
consistenly, is not a continuous image of~$\Hstar$, see~\cite{MR2425747}.
One parallel is still unresolved: is every perfectly normal continuum
a continuous image of~$\Hstar$?

For the longest time we knew of only a few (i.e., finitely many)
distinct subcontinua of~$\Hstar$;
just recently that number was raised to the maximum possible~$2^\cee$.
If $\CH$ \emph{does not} hold then this follows
from results on ultrapowers of linear orders, \cite{MR2823685},
the construction of such a family in case $\CH$ does hold
is more involved, \cite{MR3414878}.

\section{Convergence}

Another recurring theme in Alan's work is that of convergence, in particular
as related to the closure operation.
We all know that in a first-countable space a point belongs to the closure
of a set if and only if there is a sequence in that set that converges
to the point.
If one puts universal quantifiers in front of `point' and `set' then one
obtains the definition of Fr\'echet-Urysohn spaces, if one allows for
well-ordered sequences of arbitrary length then one defines radial spaces.
A weakening of the Fr\'echet-Urysohn property, sequentiality, states that
a set is closed iff
it contains all limits of converging sequences whose terms belong to the set,
allowing for well-ordered sequences of arbitrary length will define
pseudo-radial spaces.
Finally, a space is countably tight (or has countable tightness) if
$\cl A=\bigcup\{\cl B:B\in[A]^{\le\aleph_0}\}$ holds for every subset~$A$ of
the space.

The question whether compact spaces of countable tightness are sequential
is known as the Moore-Mr\'owka Problem; after Balogh showed that $\PFA$
implies a yes answer Alan established its outright consistency and
more~\cite{MR1173249}; later he showed that it is even consistent
that compact spaces of countable tightness and weight~$\aleph_1$ are
Fr\'echet-Urysohn \cite{MR1460819}.

The sequential closure $\cls A$ of a set $A$ in a topological space is equal
to the set of limits of all converging sequences with terms in~$A$
(constant sequences ensure $A\subseteq\cls A$).
In a Frechet-Urysohn space we have $\cl A=\cls A$ for all subsets; in a
sequential space one obtains $\cl A$ after iterating $\cls$ at most
$\omega_1$~many times --- the sequential order of a space is the minimum
ordinal~$\alpha$ such that all closures are reached in at most $\alpha$ steps.
Under $\CH$ there are compact sequential spaces of all possible sequential
orders; without $\CH$ much less is known: in~\cite{MR2107167}
Alan constructed an example from $\bee=\cee$ of compact space of scattered
height and sequential order~$5$ and later he showed, under~$\PFA$,
that the known constructions would not yield spaces of order
higher than~$\omega$ \cite{MR2884242}.

Efimov's problem may also be put under the heading of convergence:
it asks whether every compact Hausdorff space contains a non-trivial
converging sequence of a copy of~$\betaN$.
Fedorchuk's compact S-space from $\diamond$ is a particularly strong
counterexample: hereditarily separable and without converging sequences.
The consistency of a positive answer seems to be getting a more and more
remote possibility, not least because of examples that Alan (co-)constructed:
there is one if $\cf([\ess]^{\aleph_0},{\subset})=\ess$ and $2^\ess<2^\cee$
hold \cite{MR2182920} and also one if $\bee=\cee$ \cite{MR3164725}
(with Shelah).

The latter example helps settle another problem on convergence:
does every compact Hausdorff space contain a converging $\omega$-sequence
or a converging $\omega_1$-sequence (Hu\v{s}ek); as $\betaN$ contains
a sequence of the latter type, this question can be seen as a weakening
of Efimov's problem.
That it really is a weakening follows from results in~\cite{MR3194414}:
in various models, including the standard model for $\bee=\cee$,
compact spaces without converging $\omega_1$-sequences are first-countable.

\section{MAD families}

A look at some of the papers cited above will show that one of
Alan's tools of choice is a Maximal Almost Disjoint (MAD) family.
They are also objects of study in their own right.
A family, $\calA$, of countably infinite subsets of a set $X$ is
almost disjoint if any two distinct elements have a finite intersection
and a MAD family is maximal among such families with respect to inclusion.
One topologizes $X\cup\calA$ by making the points of $X$ isolated and by
letting the basic neighbourhoods of~$A\in\calA$ be the sets of the form
$\{A\}\cup A\setminus F$, where $F$~is finite. The resulting space,
denoted $\psi(X,\calA)$, can have many striking properties.
Alan's interest has been in the \v{C}ech-Stone remainders of such spaces,
especially in the case $\calA$~is maximal and $X=\N$.
Mr\'owka showed that this remainder can be a singleton,
and Alan and Jerry Vaughan investigated what happens for uncountable~$X$
in~\cite{MR2610447}.
One of the results that gave Alan great satisfaction deals with representable
algebras, or their duals, representable spaces: these are the zero-dimensional
compact spaces that can occur as $\psi(\N,\calA)^*$ for some MAD
family~$\calA$.
Baumgartner and Weese initiated the study of these and proved that, as
one is wont to expect, all is well under the assumption of~$\CH$:
one gets the compact zero-dimensional spaces of weight~$\cee$ or less.
In~\cite{MR1927102}, in a veritable tour de force, Alan and his student
JinYuan Zhou presented a consistent example of a mad family~$\calD$ with
two points, $a$ and $b$, in~$\psi(\N,\calD)^*$ such that the quotient of
the remainder obtained by identifying just the points $a$ and~$b$ is
not of the form $\psi(\N,\calA)^*$ for any MAD family~$\calA$.

\end{document}